\documentclass{article}
\usepackage{amsfonts}
\usepackage{amssymb}
\newcommand{\pr}{\scriptscriptstyle\perp}
\begin{document}

\begin{center}
{\Large {\bf Elementary parabolic twist
\footnote{This  work was
supported by Russian Foundation for Basic Research under the grant
N 00-01-00500. }  }}\\[5mm] { Lyakhovsky V.D.
\footnote{E-mail address: lyakhovs@pobox.spbu.ru }
 and Samsonov M.E.
\footnote{E-mail address: samsonov@heps.phys.spbu.ru}}\\[5mm]
 Theoretical Department, \\
 St. Petersburg State University,\\
198904, St. Petersburg, Russia
\end{center}

\vskip 0.5cm
\begin{abstract}
The twist deformations for simple Lie algebras $U(\mathfrak{g})$ whose
twisting elements $\mathcal{F}$ are known explicitly are usually defined on
the carrier subspace injected in the Borel subalgebra $\mathfrak{B^+(g)}$.
We solve the problem of creating the parabolic twist $\mathcal{F}_{\wp%
}$ whose carrier algebra $\mathfrak{P}$ not only covers
$\mathfrak{B^+(g)}$ but also intersects nontrivially with
$\mathfrak{B^-(g)}$. This algebra $\mathfrak{P}$ is the parabolic
subalgebra in $\mathfrak{sl(3)}$ and has the structure of the
algebra of two-dimensional motions. The parabolic twist is
explicitly constructed as a composition of the extended jordanian
twist $\mathcal{F}_{\mathcal{EJ}}$ and the new factor
$\mathcal{F}_{\mathcal{D}}$. The latter can be considered as a
special version of the j\/ordanian twist. The twisted costructure
is found for $U(\mathfrak{L})$ and the corresponding universal
$\mathcal{R}$-matrix is presented.
\end{abstract}

\section{Introduction}

There are two main sets of constant classical quasitriangular $r$-matrices
for semisimple Lie algebras $\mathfrak{g}$: nonskewsymmetric and
skewsymmetric. The first set was classified by Belavin and Drinfeld \cite{BD}%
. The corresponding $r$-matrices were associated with the subsets
$T_{1,2}$ of the root system $\Lambda(\mathfrak{g})$ and the
morphisms $\mathcal{T}$ connecting them (Belavin-Drinfeld
triples). The classification scheme for skewsymmetric $r$-matrices
was developed by Stolin \cite{STO}. In this scheme listing of
skewsymmetric solutions of the classical Yang-Baxter equation
(CYBE) was reduced to the classification of quasi-Frobenius subalgebras in $%
\mathfrak{g}$, their normalizers and the evaluation of the cohomology group $%
H^2(\mathfrak{g})$.

The existence of quantizations was proved by Drinfeld for all Lie bialgebras
\cite{EK}.
But that proof does not allow one to compute the quantization explicitly. For
a long period of time the universal $\mathcal{R}$-matrices were known only
for some classical $r$-matrices (standard or Drinfeld-Jimbo \cite{D,J}
and Cremmer-Gervais \cite{CG} for nonskewsymmetric and j\/ordanian
\cite{O},
Reshetikhin \cite{R} and GGS \cite{GGS} for skewsymmetric solutions and
some other).

It was also known that in nonskewsymmetric case some solutions of quantum
YBE can be transformed into the standard one by a special kind of twist.
The corresponding twisting element connecting Cremmer-Gervais and standard $%
\mathcal{R}$-matrices was constructed by Kulish and Mudrov \cite{KM}.
Recently Etingof, Schedler and Schiffmann \cite{ESS} had proved that such
twists exist for any quantized Belavin-Drinfeld triple and thus had solved
the problem of explicit quantization for nonskewsymmetric $r$-matrices. The
alternative expressions for the twisting elements mentioned above where
proposed in \cite{IO} by Isaev and Ogievetsky.

The quantization problem for quasitriangular Lie bialgebras with
skewsymmetric $r$-matrices was reduced by Drinfeld \cite{D83} to the solution
of the twist equations
\begin{equation}  \label{drinf}
\begin{array}{c}
\mathcal{F}_{12}(\Delta\otimes id)\mathcal{F}= \mathcal{F}%
_{23}(id\otimes\Delta)\mathcal{F}, \\
\left( \epsilon \otimes id \right) \mathcal{F} = \left( id \otimes \epsilon
\right) \mathcal{F},
\end{array}
\end{equation}
Such solutions are called twisting elements $\mathcal{F} \in U(\mathfrak{L})
\otimes U(\mathfrak{L})$ and the minimal $\mathfrak{L} \subset \mathfrak{g}$
on which  $\mathcal{F}$ is defined is
the carrier subalgebra of the twist. The well known triangular universal $%
\mathcal{R}$-matrices were found when
the corresponding solutions of  (\ref{drinf}) were constructed (j\/ordanian
 \cite{O}, Reshetikhin  \cite{R} and GGS-twists  \cite{GGS}).  Starting with
the construction of the extended j\/ordanian twists (EJT's) deeper
understanding of the peculiarities of skewsymmetric class of
solutions was achieved \cite{KLM}. The essential part of the
variety of  skew $r$-matrices (with carriers in the Borel
subalgebras $\mathfrak{B(g)} \subset \mathfrak{g}$) was explicitly
quantized. The corresponding twisting elements were found to be
the compositions ("chains") of EJT's. It becomes clear that EJT's
play the fundamental role in quantizing the skewsymmetric
$r$-matrices. On the other hand in many cases the canonical forms
of EJT's are insufficient and the precedent twisting factors of a
chain can induce deformations of the consequent ones \cite{KL}.

Let $\Lambda = \Lambda^+ \cup \Lambda^-$ be the root system of $\mathfrak{g}$
and $V_{\mathfrak{g}} = V^- \oplus  V_{\mathfrak{H}} \oplus V^+$ -- the
corresponding
triangular decomposition. While constructing the twisting elements explicitly
the most difficult is the situation when the
carrier $\mathfrak{L}$ subalgebra intersects nontrivially both $V^+$ and $V^-
$. (The intersection must be considered trivial when there exist an
automorphism  bringing $\mathfrak{L}$ into $V^+$.) In that case a chain of
twists  starting in $V^+$ deforms the costructure of the space $V_{%
\mathfrak{L}} \cap V^-$ so  that the corresponding basic elements lose their
quasiprimitivity.  This deformed costructure of $\mathfrak{L}$ must be used
in the twist equation  (\ref{drinf}) to find the twisting factors that could
enlarge the initial chain.  This leads to severe difficulties.

In this paper we consider the simplest example  of the carrier algebra $%
\mathfrak{L}$ containing the Borel subalgebra $\mathfrak{B^+(g)}$ and a
one-dimensional subspace from $\mathfrak{B^-(g)}$. This is the maximal
parabolic subalgebra $\mathfrak{P}$ in $\mathfrak{g=sl(3)}$.  As was shown
by Gerstenhaber and Giaquinto \cite{GG} the $r_{\mathfrak{P}}$-matrices for
the parabolic subalgebras lay in the  boundaries of the smooth varieties of
Cremmer-Gervais solutions $r_{CG}$ of  the modified CYBE. In particular
there exists such an element $x \in sl(3)$  that the following expression
\[
\exp(-t\mathop{\rm ad}(x))\circ r_{CG} = r_{CG} + t \cdot r_{\mathfrak{P}}
\]
is a solution of CYBE with an arbitrary parameter $t$.

Our aim is to construct the twisting element $\mathcal{F}_{\wp}$
with the carrier subalgebra
\begin{equation}
\label{p-def}
\mathfrak{P} \equiv \mathfrak{P}_1 =
 \mathfrak{B^+(sl(3))} \vdash V_{-\alpha_2} =
\mathfrak{B^+(sl(3))} \vdash V(E_{32}),
\end{equation}
where the generators $E_{ik}$ are the usual matrix units. The
problem is solved in two steps (see Section 2). First we perform
the peripheric extended twisting $\mathcal{F}_{\mathcal{P}}$
\cite{LO} with the carrier in $\mathfrak{B^+(sl(3))}$. This gives
the first two factors $\Phi _{\mathcal{J}}$ and $\Phi
_{\mathcal{E}}$ of the twisting element $\mathcal{F}_{\wp}$.
Studying the obtained costructure of the deformed algebra
$U_{\mathcal{P}}(\mathfrak{P})$ we can predict the form of the
last factor $\mathcal{F} _{\mathcal{D}}$ whose
carrier contains $E_{32}$. Our main statement is that this $%
\mathcal{F} _{\mathcal{D}}$ is a solution of the twist equations for $%
\mathfrak{P}$ with the costructure deformed by the peripheric twist $%
\mathcal{F} _{\mathcal{P}}=\Phi_{\mathcal{E}} \Phi_{\mathcal{J}}$. The
proof demonstrates the connection between the generalized
Verma's identity and the factorization of Drinfeld equation for
$\mathcal{F} _{\mathcal{D}}$. It follows that the composition
$\mathcal{F}_{\wp}= \mathcal{F}_{%
\mathcal{D}}\mathcal{F}_{\mathcal{EJ}}$ forms the twist with the parabolic
carrier subalgebra (\ref{p-def}). We call it the elementary
parabolic twist. In Section 3 the corresponding
$\mathcal{R}_{\wp}$ matrix is evaluated in the defining
representation. The applications of this result and possible
extensions of the construction are discussed.


\section{Twisting the universal enveloping algebra $U(\mathfrak{P})$}

The algebra $U(\mathfrak{sl(3)})$ will be considered in the canonical $%
\mathfrak{gl(3)}$-basis $\{E_{ij};\quad i,j=1,2,3\}$ with the relations
\[
\left[ E_{ij},E_{kl}\right] =\delta _{jk}E_{il}-\delta _{il}E_{kj}.
\]
For the Cartan subalgebra we shall use the following generators
\begin{eqnarray*}
H_{13}^{\scriptscriptstyle\perp } &\!=\!&\frac{1}{3}E_{11}-\frac{2}{3}E_{22}+%
\frac{1}{3}E_{33}, \\
H_{23}^{\scriptscriptstyle\perp } &\!=\!&\frac{2}{3}E_{11}-\frac{1}{3}E_{22}-%
\frac{1}{3}E_{33}, \\
H_{23} &\!=\!&E_{22}-E_{33}=
H_{23}^{\scriptscriptstyle\perp }-2H_{13}^{\scriptscriptstyle\perp }.
\end{eqnarray*}
Notice that the dual vectors, $H_{13}^{\scriptscriptstyle\perp
\ast }$ and $H_{23}^{\scriptscriptstyle\perp \ast }$, are
orthogonal to the roots $\alpha _{13}$ and $\alpha _{23}$
respectively.

Our aim is to construct the twist whose carrier subalgebra $\mathfrak{P}$ is
generated by the Borel subalgebra $\mathfrak{B^+(sl(3))}$ and the element $%
E_{32}$, i.e. by the parabolic subalgebra of \ $sl(3)$. The Cartan
subalgebra acts on the root generators of $\mathfrak{P}$ as follows:
\[
\begin{array}{ccccccc}
\left[ H_{13}^{\scriptscriptstyle\perp },E_{12}\right]  & \!=\! & E_{12}, &
& \left[ H_{23}^{\scriptscriptstyle\perp },E_{12}\right]  & \!=\! & E_{12},
\\[1ex]
\left[ H_{13}^{\scriptscriptstyle\perp },E_{13}\right]  & \!=\! & 0, &  &
\left[ H_{23}^{\scriptscriptstyle\perp },E_{13}\right]  & \!=\! & E_{13}, \\%
[1ex]
\left[ H_{13}^{\scriptscriptstyle\perp },E_{23}\right]  & \!=\! & -E_{23}, &
& \left[ H_{23}^{\scriptscriptstyle\perp },E_{23}\right]  & \!=\! & 0, \\%
[1ex]
\left[ H_{13}^{\scriptscriptstyle\perp },E_{32}\right]  & \!=\! & E_{32}, &
& \left[ H_{23}^{\scriptscriptstyle\perp },E_{32}\right]  & \!=\! & 0.
\end{array}
\]
For the algebra $\mathfrak{P}$ the twisting can be started by
applying the extended j\/ordanian
twist $\mathcal{F}_{\mathcal{EJ}}$ \cite{KLM} based on the Borel subalgebra $%
\mathfrak{B}_{1}$ generated by some Cartan generator $H_{1}$\ and a
generator that can be chosen to be $E_{13}$:
\begin{equation}
\label{f-ext}
\mathcal{F}_{\mathcal{EJ}}=
\exp (-E_{23} \otimes E_{12}e^{\beta \sigma_{13}} )
\exp (H_{1} \otimes \sigma_{13} )
\end{equation}
with $\sigma_{13}=\ln(1+E_{13})$ and $\left[ H_1 , E_{13}\right]
\neq 0$. Our aim is  to enlarge the canonical extended j\/ordanian
twist to incorporate the additional Cartan generator $H_{2}$\ and
the element $E_{32}$ into the carrier subalgebra $\mathfrak{P}$ of
the final twisting element. In the general case $H_{2}$\  forms
Borel subalgebras (say $\mathfrak{B}_{2}$) with both $E_{13}$ and
$E_{32}$. Evidently one can suppose that some kind of j\/ordanian
twist (based on the subalgebra $\mathfrak{B}_{2}$) might be
applicable to
the algebra $U_{\mathcal{EJ}}(\mathfrak{sl(3)})$\ deformed by $\mathcal{F}_{%
\mathcal{EJ}}$. Thus it is natural to consider the following supposition: the
twisting element with the parabolic carrier in $U(\mathfrak{sl(3)})$ has the
form
\begin{equation}
\label{f-para}
\begin{array}{lcl}
\mathcal{F}_{\wp} &=&
\mathcal{F}_{\mathcal{D}}\mathcal{F}_{\mathcal{E}}\mathcal{F}_{\mathcal{J}}
=\mathcal{F}_{\mathcal{D}}\mathcal{F}_{\mathcal{EJ}}= \\
&=&\exp (H_{2}\otimes (b\sigma_{13} +\sigma_{32}))\exp (-E_{23}\otimes
E_{12}e^{-\beta \sigma_{13} })\exp (H_{1}\otimes \sigma_{13}).
\end{array}
\end{equation}
Here $[H_{1},E_{12}]=\beta E_{12}$  and  $\sigma_{32} = \ln
(1+E_{32})$. We have two free parameters ($b$ and $\beta $) and
two yet undetermined Cartan generators: $H_{1}$ and $H_{2}$. The
ordering of tensor factors in the extension
$\mathcal{F}_{\mathcal{E}}$ isn't essential and similar
considerations work for the alternative ordering \cite{LO}.

In our notations the $r$- matrix studied in \cite{GG} looks like
\begin{equation}
\label{r-mat}
r_{\wp}=H_{23}^{\pr}\wedge E_{13}+E_{12}\wedge
E_{23}+H_{13}^{\pr}\wedge E_{32}.
\end{equation}
This $r$-matrix imposes the following conditions on the generators:
\[
\begin{array}{cc}
H_{1}+bH_{2}=H_{23}^{\scriptscriptstyle\perp }, & H_{2}=H_{13}^{%
\scriptscriptstyle\perp }.
\end{array}
\]
The factorization property $\mathcal{F}_{\wp} =
\mathcal{F}_{\mathcal{D}}\mathcal{F}_{\mathcal{EJ}}$ leads to the additional
restriction: the extended jordanian factor $
\mathcal{F}_{\mathcal{EJ}}$ in $\mathcal{F}_{\wp} $ must
survive after the generator $E_{32}$ is scaled and the scale
parameter is sent to zero. Thus we get the additional condition:
$\beta +b=1.$ Moreover the explicit form (\ref{r-mat}) of the
$r$-matrix clearly indicates that it is natural to rearrange the
factors of $\mathcal{F}_{\wp} $:
\begin{eqnarray*}
\mathcal{F}_{\wp}=
\mathcal{F}_{\mathcal{D}} \mathcal{F}_{\mathcal{E}} \mathcal{F}_{\mathcal{J}}
&=&\exp (H_{13}^{\scriptscriptstyle\perp }\otimes (b\sigma_{13}
+\sigma_{32}))\times  \\
&&\times \exp (-E_{23}\otimes E_{12}e^{-(1-b)\sigma_{13} })\exp ((H_{23}^{%
\scriptscriptstyle\perp }-bH_{13}^{\scriptscriptstyle\perp })\otimes
\sigma_{13} )\\
=\mathcal{F}_{\mathcal{D}} \mathcal{F}_{\mathcal{R}}
\mathcal{F}_{\mathcal{P}} &=&
\mathcal{F}_{\mathcal{D}} \mathcal{F}_{\mathcal{R}}
\Phi_{\mathcal{E}} \Phi_{\mathcal{J}}\\
&=&
\exp (H_{13}^{\scriptscriptstyle\perp }\otimes
(b\sigma_{13} +\sigma_{32}))  \exp
(-bH_{13}^{\scriptscriptstyle\perp }\otimes \sigma_{13} )\times
\\
&&\times\exp (-E_{23}\otimes E_{12}e^{-\sigma_{13} })\exp
(H_{23}^{\scriptscriptstyle\perp }\otimes \sigma_{13} ).
\end{eqnarray*}

We shall use this form to justify the proposed expression for the
parabolic twisting element and to fix the value of $b$. Let us
impose on $\mathcal{F}_{\wp}$ the Drinfeld condition with the
costructure of the algebra $U_{\mathcal{P}}(\mathfrak{sl(3)})$.
The latter is the result of the twist deformation performed by the
so called peripheric extended j\/ordanian twist
$\mathcal{F}_{\mathcal{P}}$ \cite{LO}:
\begin{equation}
\label{extended}
\mathcal{F}_{\mathcal{P}}=\Phi_{\mathcal{E}} \Phi_{\mathcal{J}}=
\exp (-E_{23}\otimes E_{12}e^{-\sigma_{13} })\exp
(H_{23}^{\scriptscriptstyle\perp }\otimes \sigma_{13} ).
\end{equation}
\[
\mathcal{F}_{\mathcal{P}}:U(\mathfrak{P})\longrightarrow
U_{\mathcal{P}}(\mathfrak{P}).
\]

The twisting element $\mathcal{F}_{\mathcal{P}}$ has the
4-dimensional carrier subalgebra $\mathfrak{L}\subset
\mathfrak{P}$ generated by the set $\left\{
H_{23}^{\scriptscriptstyle\perp },E_{12}, E_{23},E_{13}\right\} $.
The costructure of the twisted algebra
$U_{\mathcal{P}}(\mathfrak{P})$ is defined by the following
coproducts:
\begin{equation}
\begin{array}{lcl}
\Delta_{\mathcal{P}}(H_{13}^{\scriptscriptstyle\perp }) & \!=\! & H_{13}^{%
\scriptscriptstyle\perp }\otimes 1+1\otimes H_{13}^{\scriptscriptstyle\perp
}, \\[1ex]
\Delta_{\mathcal{P}}(E_{12}) & \!=\! & E_{12}\otimes e^{\sigma_{13} }
+e^{\sigma_{13}
}\otimes E_{12}, \\[1ex]
\Delta_{\mathcal{P}}(E_{13}) & \!=\! & E_{13}\otimes e^{\sigma_{13} }+1\otimes
E_{13}, \\[1ex]
\Delta_{\mathcal{P}}(E_{23}) & \!=\! & E_{23}\otimes e^{-\sigma_{13} }
+1\otimes E_{23}, \\[1ex]
\Delta_{\mathcal{P}}(H_{23}^{\scriptscriptstyle\perp }) & \!=\! & H_{23}^{%
\scriptscriptstyle\perp }\otimes 1+1\otimes H_{23}^{\scriptscriptstyle\perp
}+E_{23}\otimes E_{12}e^{-2\sigma_{13} }, \\[1ex]
\Delta_{\mathcal{P}}(E_{32}) & \!=\! & E_{32}\otimes 1+1\otimes
E_{32}+2H_{13}^{\scriptscriptstyle\perp }\otimes
E_{12}e^{-\sigma_{13} }.
\end{array}
\label{co-ej}
\end{equation}
Now we are ready to evaluate the constant $b$. Let us expand the element
$\mathcal{F}_{\mathcal{DR}}=\mathcal{F}_{\mathcal{D}}
\mathcal{F}_{\mathcal{R}}$ with respect to $b$:
\begin{equation}
\begin{array}{c}
\mathcal{F}_{\mathcal{DR}}=\exp (H_{13}^{\scriptscriptstyle\perp }\otimes
\sigma_{32}+\frac{b}{2}(H_{13}^{\scriptscriptstyle\perp })^{2}\otimes \lbrack
\sigma_{13} ,\sigma_{32}]+\cdots )\pmod {{b}^{2}}
\end{array}
\label{tranc}
\end{equation}
In the Drinfeld equation for $\mathcal{F}_{\mathcal{DR}}$ in the
form (\ref{tranc}) let us collect the coefficients for the first
two powers of $b$. Taking into account the relation
\[
\left[ \sigma_{13} ,E_{32}\right] =E_{12}e^{-\sigma_{13} }
\]
and the form of the last term in $\Delta_{\mathcal{P}}(E_{32})$ we
immediately come to the conclusion that $b=2$.

We want to demonstrate that the Hopf algebra
$U_{\mathcal{P}}(\mathfrak{P})$ (with the costructure
(\ref{co-ej})) can be additionally twisted by the factor
$\mathcal{F}_{\mathcal{DR}}$ that depends on $E_{32}$ and
$H_{13}^{\scriptscriptstyle\perp }$. As a result the
composition of the four factors $\Phi_{\mathcal{J}}$, $\Phi_{%
\mathcal{E}}$, $\mathcal{F}_{\mathcal{R}}$
and $\mathcal{F}_{\mathcal{D}}$ will form the twisting element
$\mathcal{F}_{\wp}=\mathcal{F}_{\mathcal{DR}}\mathcal{F}_{\mathcal{P}%
}=\mathcal{F}_{\mathcal{D}}\mathcal{F}_{\mathcal{R}}\Phi_{\mathcal{E}}
\Phi_{\mathcal{J}}$  with the carrier algebra $\mathfrak{P}$.

\textbf{Theorem}

Algebra $U_{\mathcal{P}}(\mathfrak{P})$ admits the twist with the
element
\begin{equation}
\begin{array}{lr}
\mathcal{F}_{\mathcal{DR}}=\exp (H_{13}^{\scriptscriptstyle\perp }\otimes
(2\sigma_{13} +\sigma_{32}))
\exp (-2H_{13}^{\scriptscriptstyle\perp }\otimes \sigma_{13} )&\\
&\square
\end{array}
\label{d-f}
\end{equation}

\textbf{Proof}

Consider the generalized Verma identity presented in \cite{FK} in
the following form:
\[
e^{x\ln (1+ta)}e^{(x+y)\ln (1+sb)}e^{y\ln (1+ta)}=e^{y\ln (1+sb)}e^{(x+y)\ln
(1+ta)}e^{x\ln (1+sb)}.
\]
The operators $a$ and $b$ are subject to the relations:
\[
\lbrack a,[a,b]]=0
\]
and
\[
\lbrack b,[b,a]]=0;
\]
$x,y,s,t$ are some central elements. Verma identity leads to the
property:
\[
e^{\xi \sigma_{13} }e^{\xi \sigma_{32} }e^{\xi \sigma_{13}
}=(e^{\sigma_{13} }e^{\sigma_{32} }e^{\sigma_{13} })^{\xi }.
\]
Now we need only to differentiate this expression and to put $\xi
=0$. This results in the following formula
\begin{equation}
\ln (e^{ \sigma_{13} }e^{\sigma_{32}}e^{ \sigma_{13} })= (2\sigma_{13} +
 \sigma_{32} ).  \label{verma}
\end{equation}
As a consequence the element $\mathcal{F}_{\mathcal{DR}}$ acquires
the form
\[
\mathcal{F}_{\mathcal{DR}}=\exp (H_{13}^{\scriptscriptstyle\perp }\otimes \ln
(e^{\sigma_{13} }e^{\sigma_{32}}e^{\sigma_{13} }))\exp
(-2H_{13}^{\scriptscriptstyle\perp }\otimes \sigma_{13} ).
\]
Consider the deformed coproduct $\mathcal{F}_{\mathcal{DR}}
\Delta_{\mathcal{P}}(e^{\sigma_{13}}e^{\sigma_{32}}e^{\sigma_{13} })
\mathcal{F}^{-1}_{\mathcal{DR}}$.
It is easy to verify that it is group-like. Together with the
primitivity of $H_{13}^{\scriptscriptstyle\perp }$ in
$U_{\mathcal{P}}(\mathfrak{P})$
 (\ref{co-ej}) and the special form of the last factor in
$\mathcal{F}_{\mathcal{ DR}}$ this brings us to the conclusion
that $\mathcal{F}_{\mathcal{DR}}$ satisfies the factorized
Drinfeld equations: $$ (\Delta_{\mathcal{P}}\otimes
id)({\cal{F}}_{\cal{DR}})=({\cal{F}}_{\cal{DR}})_{13}
({\cal{F}}_{\cal{DR}})_{23} $$ and $$
(id\otimes\Delta_{\cal{DR}\mathcal{P}})({\cal {F}}_{\cal {DR}})
=({\cal{F}}_{\cal{DR}})_{12} ({\cal {F}}_{\cal{DR}})_{13} $$ and
therefore satisfies also the initial equation (\ref{drinf}). This
proves the theorem. $\blacksquare$

Applying  $\mathcal{F}_{\mathcal{DR}}$ to the algebra $U_{\mathcal{P}}(%
\mathfrak{P})$ we find the final costructure of the twisted parabolic algebra:
\[
U_{\mathcal{P}}(\mathfrak{P})\stackrel{\mathcal{F}_{\mathcal{DR}}}{%
\longrightarrow }U_{\wp}(\mathfrak{P}).
\]
The coproducts for the six generators of $\mathfrak{P}$ obtain the following
form:
\begin{equation}
\begin{array}{lcl}
\Delta _{\wp}(H_{13}^{\scriptscriptstyle\perp }) & \!=\! & 1\otimes
H_{13}^{\scriptscriptstyle\perp }+(H_{13}^{\scriptscriptstyle\perp }\otimes
1)(1\otimes 1+C)^{-1}, \\[1ex]
\Delta _{\wp}(H_{23}^{\scriptscriptstyle\perp }) & \!=\! & 1\otimes
H_{23}^{\scriptscriptstyle\perp }+H_{13}^{\scriptscriptstyle\perp
}\otimes e^{-\sigma_{13} }+(E_{23}\otimes E_{12}e^{-\sigma_{13} }+
\\[1ex] &  & +((H_{23}^{\scriptscriptstyle\perp
}-H_{13}^{\scriptscriptstyle\perp })\otimes 1)(1\otimes
1+C))(1\otimes e^{\sigma_{13} }e^{\sigma_{32}})^{-1}, \\[1ex]
\Delta _{\wp}(E_{12}) & \!=\! & E_{12}
\otimes e^{\sigma_{32}}e^{\sigma_{13}
}+e^{\sigma_{13} }\otimes E_{12}+H_{13}^{\scriptscriptstyle\perp
}E_{12}\otimes E_{12}, \\[1ex]
\Delta _{\wp}(E_{{13}}) & \!=\! & (e^{\sigma_{13} }\otimes e^{\sigma_{13}
}e^{\sigma_{32}})(1\otimes 1+C)^{-1}-1\otimes 1, \\[1ex]
\Delta _{\wp}(E_{23}) & \!=\! & (E_{23}\otimes e^{-\sigma_{13} }+H_{13}^{%
\scriptscriptstyle\perp }\otimes (2H_{13}^{\scriptscriptstyle\perp }-H_{23}^{%
\scriptscriptstyle\perp })- \\[1ex]
&  & -(H_{13}^{\scriptscriptstyle\perp })^{2}\otimes
e^{-\sigma_{13} } +H_{13}^{\scriptscriptstyle\perp }\otimes
1)(1\otimes 1+C)^{-1}+ \\[1ex] &  &
+(H_{13}^{\scriptscriptstyle\perp
}(H_{13}^{\scriptscriptstyle\perp }-1)\otimes 1)(1\otimes
1+C)^{-2}+1\otimes E_{23}, \\[1ex]
\Delta _{\wp}(E_{32}) & \!=\! & E_{32}\otimes e^{\sigma_{32}}+1\otimes
E_{32}+ \\[1ex] &  &
+(E_{32}+2e^{\sigma_{32}}H_{13}^{\scriptscriptstyle\perp })\otimes
E_{12}e^{-\sigma_{13} }+ \\[1ex]
&  & +(E_{32}+e^{\sigma_{32}}H_{13}^{\scriptscriptstyle\perp })H_{13}^{%
\scriptscriptstyle\perp }\otimes (E_{12})^{2}(e^{\sigma_{13} }e^{
\sigma_{32}}e^{\sigma_{13}
})^{-1}.
\end{array}
\label{copr-par}
\end{equation}
where
\[
C=1\otimes E_{32}+H_{13}^{\perp }\otimes E_{12}e^{-\sigma_{13} }.
\]

We can consider the Hopf algebra $U_{\wp}(\mathfrak{P})$ as a
result of the integral twist deformation
\[
U(\mathfrak{P})\stackrel{\mathcal{F}_{\wp}}{\longrightarrow }U_{%
\wp}(\mathfrak{P}),
\]
where the element $\mathcal{F}_{\wp}$ can be written in the form
\[
\label{par-twist}
\begin{array}{lcl}
\mathcal{F}_{\wp}&=&
\mathcal{F}_{\mathcal{D}}\mathcal{F}_{\mathcal{EJ}}\\
& = & \exp (H_{13}^{\scriptscriptstyle\perp }\otimes (2\sigma_{13}
+\sigma_{32}))\exp (-E_{23}\otimes E_{12}e^{\sigma_{13} })\exp
(H_{23}\otimes \sigma_{13} )
\end{array}
\]
and is called the elementary parabolic twisting element.

\section{Universal element and $R$-matrix}

The parabolic twist $\mathcal{F}_{\wp}$ can be supplied with two
natural parameters corresponding to two j\/ordanian-like
deformations. The
following rescaling of the generators is an automorphism of the algebra $%
\mathfrak{P}$:
\[
E_{13}\longrightarrow \xi E_{13}, \quad E_{32}\longrightarrow \zeta E_{32},
\quad E_{12}\longrightarrow \xi \zeta E_{12}, \quad E_{23}\longrightarrow
\frac{1}{\zeta} E_{23}.
\]
It induces the parametrization of the twisting element
\begin{equation}
  \label{par-param}
\begin{array}{l}
\mathcal{F}_{\wp}(\xi, \zeta)=\\
\exp(H_{13}^{\scriptscriptstyle%
\perp}\otimes (2\sigma_{13}(\xi)+\sigma_{32}(\zeta))) \exp(- \xi E_{23}\otimes
E_{12}e^{\sigma_{13}(\xi)}) \exp(H_{23}\otimes\sigma_{13}(\xi)),
\end{array}
\end{equation}
where $\sigma_{ij}(\xi)=\ln(1+ \xi E_{ij})$.

The result of the parabolic twisting with the element
(\ref{par-param}) is the
2-dimensional smooth variety of Hopf algebras $U_{\wp}(\mathfrak{P}%
;\xi,\zeta)$. The parameters are independent. In the limit points
we get the ordinary twists:
\[
\begin{array}{c}
\mathcal{F}_{\wp}(\xi, \zeta) \stackrel{\zeta \rightarrow 0}{%
\longrightarrow} \mathcal{F}_{\mathcal{P}}(\xi), \\
\mathcal{F}_{\wp}(\xi, \zeta) \stackrel{\xi \rightarrow 0}{%
\longrightarrow} \mathcal{F}_{\mathcal{J}}(\zeta).
\end{array}
\]
Here the first limit is the parametrized peripheric extended twist (\ref
{extended}), the second is the j\/ordanian twist for the Borel subalgebra
generated by $\left\{H_{13}^{\scriptscriptstyle\perp}, E_{32}\right\}$ with
the twisting element $\mathcal{F}_{\mathcal{J}}(\zeta)= e^{H_{13}^{%
\scriptscriptstyle\perp} \otimes \sigma_{32}(\zeta )}$.

The universal $\mathcal{R}$-matrix for $U_{\wp}(\mathfrak{P}; \xi,
\zeta)$ is defined by the standard expression:
\begin{equation}  \label{rrmat-param}
\begin{array}{c}
\mathcal{R}_{\wp}(\xi, \zeta) =\left(\mathcal{F}_{\wp}(\xi, \zeta
)\right)_{21} \left(\mathcal{F}_{\wp}(\xi, \zeta )
\right)^{-1}= \\
\exp((2\sigma_{13}(\xi)
+\sigma_{32}( \zeta))\otimes H_{13}^{\scriptscriptstyle%
\perp}) \exp(-\xi E_{12}e^{\sigma_{13}(\xi)}\otimes E_{23})
\exp(\sigma_{13}(\xi)\otimes H_{23})\times \\
\times \exp(-H_{23}\otimes\sigma_{13}(\xi)) \exp(\xi E_{23}\otimes
E_{12}e^{\sigma_{13}(\xi)}) \exp(-H_{13}^{\scriptscriptstyle\perp}\otimes(2%
\sigma_{13}(\xi)+\sigma_{32}( \zeta)).
\end{array}
\end{equation}
If we choose the parameters to be proportional ($\zeta = \eta \xi
$) then the expression (\ref{rrmat-param}) can be considered as a
quantized version of the classical $r$-matrix
\[
r _{\wp}(\eta) = H_{23}^{\scriptscriptstyle\perp} \land E_{13} +
E_{12} \land E_{23} + \eta H_{13}^{\scriptscriptstyle\perp}
\land E_{32}.
\]
In the fundamental representation the $R$-matrix has the form:
\[
\begin{array}{c}
R_{\wp}=1\otimes 1+ \\[1ex] (E_{13}\land (\frac 23 E_{11}-\frac 13
E_{22}-\frac 13 E_{33}) + E_{23}\land E_{12})\xi +
\\ +\frac 29 (E_{13}\otimes E_{13}){\xi}^{2}+ \\[1ex]
+ (E_{32}\land (\frac 13 E_{11}-\frac 23 E_{22}+\frac 13 E_{33}))\zeta + \\%
[1ex]
+ \frac 29 (E_{32}\otimes E_{32}){\zeta}^{2}+ \\[1ex]
+\frac 13 (E_{12}\otimes (\frac 13 E_{11}-\frac 23 E_{22} +\frac 13 E_{33})+
(\frac 13 E_{11}-\frac 23 E_{22}+\frac 13 E_{33})\otimes E_{12} + \\[1ex]
+ \frac 13 (E_{13}\otimes E_{32}+E_{32}\otimes E_{13}))\zeta\xi - \\[1ex]
-\frac {2}{81} (E_{12}\otimes E_{12}){\zeta}^{2}{\xi}^{2}+\frac {2}{27}
(E_{12}\land E_{32})\xi {\zeta}^{2} + \\[1ex]
+ \frac {1}{27}(E_{12}\land E_{13}){\xi}^{2}\zeta.
\end{array}
\]

\section{Conclusions}

We have demonstrated that with the help of the peripheric EJT's
and the j\/ordanian-like factors the carrier subalgebra can be
enlarged so that the negative sector of the Cartan decomposition
is partially occupied. The elementary parabolic twist
(\ref{par-twist}) can be applied to any Lie algebra that contains
the parabolic carrier subalgebra $\mathfrak{P}$. The algebras
$\mathfrak{sl(3)}$ and $\mathfrak{G}_2$ are the only ones with
this property among the simple Lie algebras with
$\mathrm{rank}(\mathfrak{g})=2$. Any simple algebra whose rank is
greater than 2 contains $\mathfrak{P}$ and consequently can be
twisted by $\mathcal{F}_{\wp}$. There are substantial reasons to
suppose that constructions similar to the composition
$\mathcal{F}_{\wp}=
\mathcal{F}_{\mathcal{D}}\mathcal{F}_{\mathcal{EJ}}=
\mathcal{F}_{\mathcal{D}}\mathcal{F}_{\mathcal{R}}\mathcal{F}_{\mathcal{P}}$
exist for the parabolic subalgebras of larger dimensions. So the
elementary parabolic twist can be considered as the first step
towards the invasion of the negative sector of the Cartan
decomposition.

There exist other twists with 6-dimensional carriers, for example
the one constructed in \cite{KL}. In the latter case the twist can
be presented as a chain containing the extended j\/ordanian twist
and the j\/ordanian twist on the deformed carrier space. It is
important to know whether the similar interpretation can be used
in the case of the parabolic twist.

Any maximal parabolic subalgebra $\mathfrak{P}_i$ with the missing
negative simple root $\lambda_i$ can be considered as a semidirect
product of the maximal simple subalgebra in $\mathfrak{P}_i$ and
the ideal generated by the basic elements $E_{\alpha_k}$ whose
positive roots $\alpha_k$ contain the simple component
$\lambda_i$. This means that the corresponding Hopf algebras
$U_{\wp_i}(\mathfrak{P}_i)$ are the examples of the algebras of
motion $\mathfrak{g}$ quantized by twists whose carriers
coincide with $\mathfrak{g}$. In the case considered here the algebra $U(%
\mathfrak{P})$ is the universal enveloping of the semidirect product $%
\mathfrak{P \approx gl(2) \vdash t(2)}$, that is the algebra of
two-dimensional motions.

\end{document}